%% file: main.tex
\documentclass[accept,onefignum,onetabnum]{siamart190516}
\UseRawInputEncoding

\input{ex_shared}

\ifpdf
\hypersetup{
  pdftitle={A Data Analysis Study On Human Liver Blood Circulation},
  pdfauthor={T. Wu, K. Liu, E. Fainman, and Y. Wang}
}
\fi

\begin{document}

\maketitle

\begin{abstract}
  The liver has a unique blood supply system and plays an important role in the human blood circulatory system. Thus, hemodynamic problems related to the liver serve as an important part in clinical diagnosis and treatment. Although estimating parameters in these hemodynamic models is essential to the study of liver models, due to the limitations of medical measurement methods and constraints of ethics on clinical studies, it is impossible to directly measure the parameters of blood vessels in livers.
  Furthermore, as an important part of the systemic blood circulation, livers' studies are supposed to be in conjunction with other blood vessels. In this article, we present an innovative method to fix parameters of an individual liver in a human blood circulation using non-invasive clinical measurements. The method consists of  a 1-D blood flow model of human arteries and veins, a 0-D model reflecting the peripheral resistance of capillaries and a lumped parameter circuit model for human livers. We apply the finite element method in fluid mechanics of these models to a numerical study, based on non-invasive blood related measures of 33 individuals. The estimated results of human blood vessel characteristic and liver model parameters are verified from the perspective of Stroke Value Variation, which shows the effectiveness of our estimation method.
\end{abstract}

\begin{keywords}
  1-D model of blood flow, liver model, Finite Element Method, data analysis
\end{keywords}

\begin{AMS}
\end{AMS}

\section{Introduction}
For vertebrates, the liver is an important metabolic organ---it receives 20\% of the amount of bleeding from the heart pump in each cardiac cycle, which significantly affects the blood circulation of the whole body. Moreover, the blood supply of the liver is very special. The hepatic portal vein (PV) and hepatic artery (HA) are the two main blood sources of the liver\cite{Garcea2009}. The blood provided by the PV accounts for about 70-80\% of the total blood inflows to the liver, and is rich in nutrients required for liver metabolism; while the HA can control and regulate the blood flow through the liver\cite{Bonfiglio2010} . According to the structure of blood vessels in the liver, we can divide the liver into several liver lobes anatomically. The hepatic PV and HA differentiate continuously in the liver lobe, forming a dense blood vessel branch network, and complete the material exchange with liver cells in the hepatic sinusoids, and then flow into the peripheral hepatic venules. The peripheral hepatic venules converge to form the hepatic vein, and the blood in the liver flows out of the liver through the hepatic vein and enters the inferior vena cava.

The structure of the hepatic blood circulation system is complex, and hence the relevant hemodynamic problems are essential in clinical diagnosis and treatment. For example, in partial hepatectomy, in order to accurately control the resection ratio to ensure that the remaining liver after the operation can maintain the normal metabolic activities of the human body, it is necessary to study the relationship among the liver resection ratio, the blood pressure and blood flow of the blood vessels in the liver \cite{Audebert2017} . In addition, the hemodynamic changes of the liver blood supply relate to the development and deterioration of several liver diseases\cite{Berzigotti2004}. Therefore, the hemodynamics of the liver has always been a concern of scientists, physicians and surgeons.
Furthermore, as an important organ in the human body, the liver involves a complex blood vessel network in conjecture with the blood circulation system in human body.
Based on the human blood circulation system, blood flow starts from the inferior vena cava, flows through the heart, then flows through the liver, and finally flows through the three leaf tissues in the liver into the inferior vena cava. Therefore, in order to study the liver blood supply model, it is necessary to fully consider the relationship among the various organs and tissue structures of the human body, when establishing a model for the blood circulation of the human body.
Unfortunately, it is difficult to directly measure the characteristics of blood vessel passing through a liver using non-invasive technical means such as electromagnetic probes and Doppler ultrasound, especially when the human body is performing normal physiological activities.
However, when using blood vessel models and liver models for research, the setting of model parameters is very important.
Therefore, existing studies generally used the liver blood supply research data of pigs, mice and other animals to compare the blood supply status of the human liver. Due to the differences in the physiological structure between animals and humans, the conclusions drawn from animal experimental studies are often different from the real situation of the human body.
In addition, the individual differences of different human bodies are ignored, when setting relevant parameters based on the results of animal experiments. This can cause errors in clinical studies and is not reliable for clinical practice.

To bridge the gap, we focus on the estimation methods of 1-D blood vessel model and liver model based on real clinical trial data.
Specifically, we use mathematical and physical methods to establish a blood vessel model suitable for non-invasive measures of humans, simulate the blood flow in the human body through the finite element method (FEM), and estimate the blood flow and blood pressure status of the liver when the human body under normal physiological conditions.
In other words, based on the acquired non-invasive clinical measures of 33 individuals in healthy conditions, we use a 1-D vascular model to model the blood vessels, and design related FEMs to simulate the blood flow status of abdominal aorta, radial artery, hepatic artery and other blood vessels; followed by estimating vascular characteristics to be used to numerically study the human liver, and obtaining the estimated values of the liver model parameters. Finally, we verify the estimated results with the observed measures in our clinical data.
Our study is designed typically to be based on the non-invasive measures of human bloods, which provides comprehensive research foundations and decision-making basis for clinical practice to increase the quality of liver-related operations and reduce the consumption of medical equipment and expenses.

The rest of the paper is organized as follows. We summarize relevant studies in \cref{sec:conclusions}, and then present blood vessel models with boundary conditions and liver models in \cref{sec:model}, where the FEMs of each model are also discussed.
Our algorithms based on FEMs is in \cref{sec:alg} with numerically experimental results, and the conclusions follow in
\cref{sec:conclusions}.

\section{Literature Review}
\label{sec:lit_review}

Our study relates to the vascular hemodynamic models and the theories of the liver model, and the FEM in fluid mechanics. In this section, we discuss these relevant existing research  from the three aforementioned aspects.

The vascular hemodynamic model can be traced back to 1775, when Euler first proposed a one-dimensional (1-D) model of the human arterial system, mainly based on the conservation of mass and momentum of inviscid fluids. At the same time, he also proposed an equation reflecting the relationship between the pressure at any point in the blood vessel and the cross-section to improve the model, but Euler himself and related researchers did not find a suitable solution for the model\cite{Euler1763}.
In 1895, Frank proposed the famous Windkessel (elastic cavity) theory, that is, the arterial system is regarded as an elastic cavity, and the microcirculation involved in the small arteries is considered as the peripheral resistance of the elastic cavity\cite{Westerhof2009}. The theory of elastic cavity is an important progress in the research of cardiovascular system model, which has a wide range of applications and has been continuously developed and improved.
With the development of partial differential equations, Riemann proposed the characteristic line method \cite{Riemann1892} suitable for general equations. Although Euler's 1-D model of blood vessel equations is a nonlinear system similar to shallow water equations, under the physiological conditions of the human arterial system, the nonlinearity of the 1-D model equations is very weak. Whomersley linearizes the equations and combines them to obtain the traveling wave solution of the blood flow equation \cite{Womersley1957}. Later on,  Stettler et al. (1981) used the characteristic variables of the 1-D blood vessel model to study the phenomenon of arterial blood flow propagation for the first time \cite{Stettler1981}.
Recently, with the rapid advancement of medical imaging technology, researchers can create three-dimensional (3-D) anatomical models to simulate blood flow from 3-D medical imaging data.
For example, due to the beating of the heart, the arterial blood flow is affected by the pulse wave, and the 3-D blood vessel model can capture the pulse wave phenomenon very well \cite{Janela2010}.
To reflect the deformation of blood vessels under the action of blood pressure, Colciago et al.(2014) built an fluid-structure interaction model \cite{Colciago2014}.
Because the numerical simulation of complex three-dimensional structural fluid problems involves huge computational costs, in actual applications, there are a very limited amount of 3-D models to numerical simulate the blood circulation system, such as the carotid artery bifurcation.
In the meanwhile, problems of the boundary conditions in the FEM numerical also arise. George Papadakis' research shows that the outflow boundary conditions have a significant impact on the pressure, blood flow velocity and other characteristics \cite{Papadakis2009}.
To solve the problem of boundary conditions, researchers introduced a 0-D total parameter model to connect the 3-D or 1-D model parts in the blood circulation system. This spatial multi-scale modeling approach uses 3-D or 1-D models to describe in details the dynamics of a specific local area.
For instance,  Formaggia et al. (1999) studied the coupling of the 1-D/0-D model, where they used 0-D model to describe internal organs, capillaries and other parts \cite{Formaggia1999}. Lagana et al. studied the coupling of 3-D and 0-D models \cite{Lagana2002}.
More recently, Quarteroni et al. (2003) proved the he existence and uniqueness of the solution to these multi-scale models \cite{Quarteroni2003}.
In most cases, the 1-D model can fully reflect the characteristics of arterial and venous blood flow, and the computational cost of numerical simulation is much smaller than that of the 3-D model. Therefore, we use a 1-D blood vessel model to simulate the propagation of blood flow in the blood vessel, and strike a balance between accuracy and computational cost.

The main feature of the blood circulatory system of the liver is that there are two blood flow channels, one is the PV that provides deoxygenated blood to the liver, and the other is the HA that provides oxygen-rich blood to the liver.
Different researchers have conducted related studies on liver models on different scales.
Chu and Reddy (1992) proposed a lumped model of visceral and liver blood circulation to illustrate the relationship among hepatic vein pressure increase, vasoconstriction and hepatic interstitial fluids\cite{Chu1992}.
Van Der Plaats (2004) used an organ-scale liver model to study cryogenic machine perfusion\cite{VanDerPlaats2004}.
3-D computational fluid dynamics  simulation was used in the portal vein blood flow simulation related to hepatectomy\cite{Ho2012}. Childress and Kleinstreuer (2014) completed 3-D computational fluid dynamics  simulation of hepatic artery blood flow for flexible and rigid blood vessels to study drug targeting\cite{Childress2014}.
Moreover, Ricken et al. (2015) proposed a porous liver model to study glucose transport and metabolism\cite{Ricken2015}. Other researchers also used this model to study liver cirrhosis\cite{Peeters2015}, the effect of liver deformation on blood pressure and blood flow\cite{Bonfiglio2010}.
Recently, some researchers proposed  a 0-D lumped parameter model to simulate the liver and couple it with blood vessel models in other parts. Audebert et al. (2017) found that a parallel circuit model can be used to represent the liver lobe. The liver can be simulated as a 0-D lumped parameter model composed of circuit models and coupled with a 1-D blood vessel model to form a complete systemic blood flow model. Based on the experimental data of pig livers, they successfully simulated the changes of hepatic artery blood flow before and after hepatectomy \cite{Audebert2017partial}.
To reflect the connection and mutual influence between the blood flow of the liver and the blood vessels in other parts of the human body, we apply the parallel circuit model proposed by Audebert et al (2017) \cite{Audebert2017partial}.

The fluid mechanics equations are mainly based on the Navier-Stokes equation, the principle of which is the basic conservation law in physics. Because the analytical solutions of Navier-Stokes equations cannot be obtained, researchers generally use computational fluid dynamics tools to study the numerical solutions of differential equations.
When studying the 1-D blood vessel model, it is necessary to use appropriate numerical calculation methods to solve the equations. Researchers have proposed many different numerical methods to study fluid dynamics models in biomedicine, including the characteristic method\cite{Schaaf1972}, the discontinuous Galerkin method\cite{Sherwin2003}, and the Taylor-Galerkin method derived from the Galerkin method\cite{Formaggia2003}, finite volume method \cite{Toro2013} and so on.
Numerical methods such as Godunov are used in venous blood vessels\cite{Brook2002}.
Boileau et al. (2015) compared the performance of six different numerical methods in the solution of the arterial vessel model\cite{Boileau2015}.
To ensure the non-negativity of the cross-sectional area of the blood vessel in the numerical solution process, Audusse et al. (2000) proposed a dynamic solution method, which is suitable for both arterial vessels and venous collapsible tubes\cite{Audusse2000}.
In this paper, we use the dynamic solution method for Saint-Venant problems ,proposed by Audusse et al.(2000) \cite{Audusse2000}, to discretize the variables of the 1-D blood vessel model and solve the model numerically.

\section{Models}\label{sec:model}

We present notations in \cref{tab:notation}.

\begin{table}[tbhp]
{\footnotesize
  \caption{Table of notations}\label{tab:notation}
\begin{center}
  \begin{tabular}{cc|cc} \hline
  \multicolumn{4}{l}{1-D vessel model} \\ \hline
   x      & position             & l       & length of the vessel segment     \\
    t      & time             & p(x,t)  & blood pressure         \\
    A(x,t) & area of the vessel cross section     & $P_0$   & reference blood pressure    \\
    $A_0$  & reference area of a cross section & $\beta$ & Vascular elasticity parameters \\
    $Q(x,t)$ & blood flow           & $\rho$  & blood density    \\
    $v(x,t)$ & blood flow velocity         & $K_f$   & friction coefficient     \\
    $\mu$  & blood viscosity coefficient     &  $V(x,t)$       &     blood volume         \\ \hline
    \multicolumn{4}{l}{Liver model} \\ \hline
   $P_a$ & hepatic artery (HA) blood pressure & $P_v$ & hepatic vein (HV) blood pressure \\
        $Q_{a,i}$ & HA blood flow at liver lobe i & $Q_{v,i}$ & HV blood flow at liver lobe $i$\\

        $R_{ha}^i$ & HA resistance value at liver lobe $i$ & $R_l^i$ & HV resistance value at liver lobe $i$ \\
        $P_{pv}$ & portal vein (PV) blood pressure & $P_{t,i}$ & blood pressure of liver lobe $i$  \\

        $Q_{pv,i}$ & PV blood flow at liver lobe $i$ & $C_l^i$ & capacitance in liver lobe $i$ circuit \\

        $R_{pv}^i$ & PV resistance value at liver lobe i & $M_i$ & weight of liver lobe $i$ \\
        $Q_{a}$ & total blood flow of the HA tree & $Q_{v}$ & total blood flow of the HV tree \\
        $Q_{pv}$ & PV tree blood flow & $R_{ha}$ & total resistance value of the HA tree \\
        $R_l$ & total resistance value of the HV tree & $R_{pv}$ & total resistance value of the HV tree \\
        $P_{t}$ & liver blood pressure & $M$ & total liver mass \\
        $C_l$ & total circuit capacitance & & \\\hline
  \end{tabular}
\end{center}
}
\end{table}

\subsection{Modeling blood circulation in vessels}

\begin{figure}[htbp]
  \centering
\label{fig:vessel_1D}
\includegraphics[scale = 0.15]{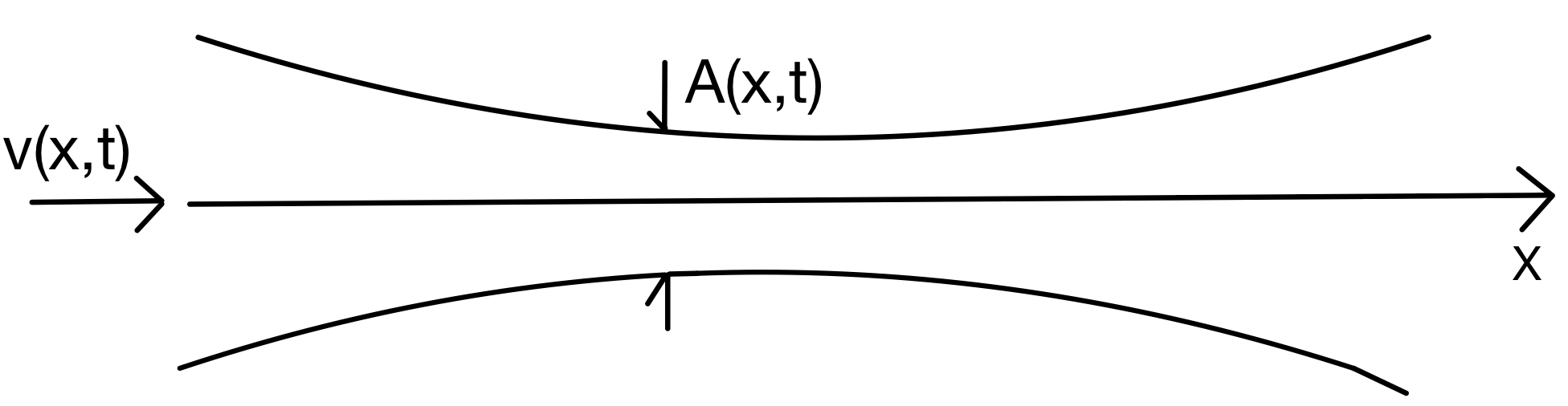}
  \caption{Illustrative model on the blood circulation in a vessel segment.}
\end{figure}

We use 1-D blood vessel model \cref{fig:vessel_1D} to show blood circulation in a segment of arteries or veins.

\subsubsection{Hemodynamic equations}\label{sec:vessel_model}
We assume blood flow is in-compressible in vessels, hence the blood density $\rho$ and viscosity $\mu$ are constant. In the meantime, in non-capillary blood vessels, the flow of blood follows the classic Newtonian fluid mechanics theorem, and there is no blood leakage. According to the law of conservation of mass, in a section of blood vessel, the change of blood mass per unit time is equal to the mass of net outflow. That is, we have the following equation \cite{Sherwin2003}
\begin{equation*}
  \rho \frac{\mathrm{d}V(t)}{\mathrm{d}t} + \rho Q(l,t) - \rho Q(0,t) = 0.
\end{equation*}
We write the integral form of the blood volume $V(t)$ and blood flow $Q(l,t) - Q(0,t)$ and derive
\begin{equation}\label{eq:1d_blood_mass}
  \frac{\partial A}{\partial t} + \frac{\partial Q}{\partial x} = \frac{\partial A}{\partial t} + \frac{\partial vA}{\partial x} = 0
\end{equation}
According to the momentum equation, the sum of the change rate of momentum in a section of blood vessel and the net flux of momentum outside the section of momentum is equal to the external force acting on this section of blood vessel; that is,
\begin{equation}
  \frac{\mathrm{d}}{\mathrm{d}t} \int_0^l \rho Q \mathrm{d}x + (\alpha \rho Qv)_l -(\alpha \rho Qv)_0 = F
\end{equation}
Let $F$ be the force in the x-axis direction of the blood vessel, and $\alpha$ the correction coefficient of momentum flux. Without loss of generality, we set $\alpha$ as 1. Then we can derive the momentum equation of the 1-D blood vessel model as
\begin{equation}\label{eq:1d_blood_momentum}
  \frac{\partial Q}{\partial t} + \frac{\partial(\alpha Qv)}{\partial x} = -\frac{A}{\rho}\frac{\partial p}{\partial x} - f + Ag
\end{equation}
In the friction function $f = K_f v(x,t)$, where coefficient of friction $K_f$ is a constant.
Because \cref{eq:1d_blood_mass} and \cref{eq:1d_blood_momentum} only reflect the law of conservation of mass and momentum, respectively, yet do not consider the expansion of blood vessels. Therefore, we introduce the relationship equation between the average blood pressure and the cross-sectional area of the blood vessel:
\begin{equation}\label{eq:1d_blood_relation}
  p(x,t) = P_0(x) + \phi (A(x,t),A_0(x),\beta (x)),
\end{equation}
where $A_0(x)$ is the reference area of cross section, $P_0(x)$ is the reference blood pressure at the reference point, and $\phi (A(x,t),A_0(x),\beta (x))$ are known parameters in the tube law function. Now we can define the tube law function for a section of artery
\begin{equation*}
  \phi (A,A_0,\beta_a) = \beta_a(A^{0.5}-A_0^{0.5}).
\end{equation*}
For a section of vein, the tube law function is
\begin{equation*}
  \phi (A,A_0,\beta_v) = \beta_v(\frac{A}{A_0}^{10}-\frac{A}{A_0}^{-1.5}).
\end{equation*}
With \cref{eq:1d_blood_mass}, \cref{eq:1d_blood_momentum} and \cref{eq:1d_blood_relation}, we obtain the 1-D blood vessel model equation set:
\begin{equation}\label{eq:1d_blood_set}
  \left\{
  \begin{aligned}
     & \frac{\partial A}{\partial t} + \frac{\partial Q}{\partial x} = 0                                                          \\
     & \frac{\partial Q}{\partial t} + \frac{\partial(\alpha Qv)}{\partial x} +\frac{A}{\rho}\frac{\partial p}{\partial x} = Ag-f \\
     & p(x,t) = P_0(x) + \phi (A(x,t),A_0(x),\beta (x))
  \end{aligned}
  \right.
\end{equation}

\subsubsection{FEM}\label{ch:FEM_1D_vessel}
Among a variety of numerical solutions for 1-D vascular equations in related studies, we consider a dynamic scheme that is suitable for both venous and arterial systems. By introducing a linear and microscopic dynamic model equivalent to the macroscopic partial differential equations \cref{eq:1d_blood_set}, we obtain the solution to the 1-D blood vessel model equations. We present the derivation process of the finite element solution method as follows\cite{Audusse2000}.

Let  $\chi(\omega) = \frac{1}{2\sqrt{3}}\mathbf{1}_{|\omega|\leq \sqrt{3}}$ be a real function with compact support on $\mathbb{R}$ and satisfy the following properties:
\begin{equation*}
  \left\{
  \begin{aligned}
     & \chi (-\omega) =   \chi (\omega) \geq 0                                                                      \\
     & \int_{\mathbb{R}}\chi (\omega) \mathrm{d}\omega = \int_{\mathbb{R}}\omega^2\chi (\omega) \mathrm{d}\omega =1
  \end{aligned}
  \right.
\end{equation*}
Suppose a distribution function
\begin{equation*}
  M(x,t,\xi) = \frac{A}{\gamma}\chi(\frac{\xi - v}{\gamma})
\end{equation*}
where
$\gamma^2 = \frac{1}{A}\int_{\epsilon A_0}^{A(x,t)}c^2(a)\mathrm{d}a$, and $c$ is the wave velocity of blood flow. For arteries, the wave velocity of blood flow propagation is defined as:
\begin{equation}\label{eq:1D_blood_velocity} 
  c^2 = \frac{\beta_a}{2\rho}\sqrt{A(x,t)}
\end{equation}
The wave velocity of blood flow in venous vessels can be referred to related literature \cite{Suzuki2019}. The relationship between between the distribution function $M$ and the solution to macro-dynamic equations \cref{eq:1d_blood_set} $(A,v)$ is
\begin{equation}\label{eq:1d_blood_FEM_M_rhs0} 
  \int_\mathbb{R} M \mathrm{d}\xi = A,\qquad
  \int_\mathbb{R}\xi M \mathrm{d}\xi = Av
\end{equation}

When the mechanical source term at the right end of the momentum equation of \cref{eq:1d_blood_set} is zero, we can prove that $(A,v)$ is the solution to \cref{eq:1d_blood_set} if and only if $M$ is the solution to the following equation
\begin{equation*}
  \frac{\partial M}{\partial t} + \xi\frac{\partial M}{\partial x} = \mathcal{G}(x,t,\xi)
\end{equation*}
where function $\mathcal{G}(x,t,\xi)$ satisfies
\begin{equation*}
  \int_\mathbb{R} \mathcal{G}(x,t,\xi)\mathrm{d}\xi = \int_\mathbb{R} \xi \mathcal{G}(x,t,\xi)\mathrm{d}\xi = 0
\end{equation*}
Let $\Delta x$ and $\Delta t$ be the space and time steps, $x_i = i\Delta x$, $t_n = n\Delta t$. Let $(A_i^n,v_i^n)$ be the numerical solution of $(A(x_i,t_n),v(x_i,t_n))$. We denote the discrete form of function $ M(x,t,\xi)$ as $M_i^n$ and $\gamma$ as $\gamma_i^n$. Specifically,
\begin{equation*}
  M_i^n(\xi) = \frac{A_i^n}{\gamma_i^n}\chi (\frac{\xi - v_i^n}{\gamma_i^n}),\qquad \gamma_i^n = (\frac{1}{A_i^n}\int_{\xi A_0}^{A_i^n}c^2(a)\mathrm{d}a)^{0.5}
\end{equation*}
Suppose $M_{i+\frac{1}{2}}^n = M_i^n\mathbf{1}_{\xi \geq 0} + M_{i+1}^n\mathbf{1}_{\xi \leq 0}$, the recursive expression of the discrete form of the distribution function $M$ becomes
\begin{equation}\label{eq:1d_blood_FEM_M_recurse} 
  M_i^{n+1} = M_i^n - \frac{\Delta t}{\Delta x}\xi(M_{i+\frac{1}{2}}^n -M_{i-\frac{1}{2}}^n )
\end{equation}
From \cref{eq:1d_blood_FEM_M_rhs0} and \cref{eq:1d_blood_FEM_M_recurse}, we have
\begin{equation*}
  X_i^{n+1} = \left(
  \begin{aligned}
       & A_i^{n+1}          \\
       & A_i^{n+1}v_i^{n+1} \\
    \end{aligned}
  \right) = \int_{\mathcal{R}} \left(
  \begin{aligned}
      1   \\
      \xi \\
    \end{aligned}
  \right) M_i^{n+1}\mathrm{d}\xi
\end{equation*}
After differentiating $X$ with respect to time $t$, we obtain the following recurrence formula
\begin{equation}\label{eq:1d_blood_FEM_Xrecurse} 
  X_i^{n+1} = X_i^{n} - \frac{\Delta t}{\Delta x}\int_\mathbb{R}\xi\left(
  \begin{aligned}
    1   \\
    \xi \\
  \end{aligned}
  \right)(M_{i+\frac{1}{2}}^n-M_{i-\frac{1}{2}}^n)\mathrm{d}\xi
\end{equation}
According to the definition of the real function $\chi(\omega)$, we have the integral of the right hand side of \cref{eq:1d_blood_FEM_Xrecurse}
\begin{equation}\label{eq:define_UV}
\begin{aligned}
 & U_{q,i}^n =  \int_{\xi\geq 0}[\xi^q\frac{A_i^n}{\gamma_i^n}\chi(\frac{\xi - v_i^n}{\gamma_i^n})]\mathrm{d}\xi = \frac{A_i^n}{2\sqrt{3}\gamma_i^n (q+1)}[(\xi\gamma_i^n + v_i^n)^{q+1}]_{\xi = max(\frac{-v_i^n}{\gamma_i^n},-\sqrt{3})}^{\xi = max(\frac{-v_i^n}{\gamma_i^n},\sqrt{3})},\\
 & V_{q,i}^n =  \int_{\xi\leq 0}[\xi^q\frac{A_i^n}{\gamma_i^n}\chi(\frac{\xi - v_i^n}{\gamma_i^n})]\mathrm{d}\xi = \frac{A_i^n}{2\sqrt{3}\gamma_i^n (q+1)}[(\xi\gamma_i^n + v_i^n)^{q+1}]_{\xi = min(\frac{-v_i^n}{\gamma_i^n},-\sqrt{3})}^{\xi = min(\frac{-v_i^n}{\gamma_i^n},\sqrt{3})},
  \end{aligned}
\end{equation}
where $q = 1,2$.
With \cref{eq:define_UV}, we can write update blood parameters using FEM as
\begin{equation}\label{alg:update_AXv}
\begin{aligned}
& A_{i}^{n+1} = A_{i}^n - \frac{\Delta t}{\Delta x}[(U_{1,i}^n + V_{1, i+1}^n) - (U_{1,i-1}^n + V_{1,i}^n)],\\
&X_{i}^{n+1} = X_{i}^{n} - \frac{\Delta t}{\Delta x}[(U_{2,i}^n + V_{2,i+1}^n) - (U_{2,i-1}^n + V_{2,i}^n)]+\Delta t(gA_{i}^n - K_fv_{i}^n),\\
& v_{i}^{n+1} = X_{i}^{n+1}/A_{i}^{n+1}.
\end{aligned}
\end{equation}
Moreover, the FEM scheme for the blood pressure \cref{eq:1d_blood_relation} can be written as
\begin{equation}\label{eq:alg_P_update}
    P^n = P_0 + \beta(\sqrt{A_{L-1}}-\sqrt{A_0})
\end{equation}

When the mechanical source term at the right hand side of the momentum equation is not zero, we adopt the following recursive form of the FEM solution:
\begin{equation*}
  X_i^{n+1} = X_i^{n} - \frac{\Delta t}{\Delta x}\int_\mathbb{R}\xi\left(
  \begin{aligned}
    1   \\
    \xi \\
  \end{aligned}
  \right)(M_{i+\frac{1}{2}}^n-M_{i-\frac{1}{2}}^n)\mathrm{d}\xi + \Delta t\left(
  \begin{aligned}
       & 0                           \\
       & gA_i^n - f(A_i^n,A_0,v_i^n) \\
    \end{aligned}
  \right)
\end{equation*}
The stability condition of FEM requires $\Delta t \cdot max(|v_i^n| + \omega_0\gamma_i^n)\leq \Delta x$, which is valid as long as the cross section area $A$ is non-negative when initial value $t=0$, this guarantees $A_i^n\geq 0, \forall i,n$.

\subsection{Modeling blood flow at vascular bifurcation}\label{ch:vessel_bif_model}

\begin{figure}[htbp]
  \centering
\label{fig:vessel_bif}
\includegraphics{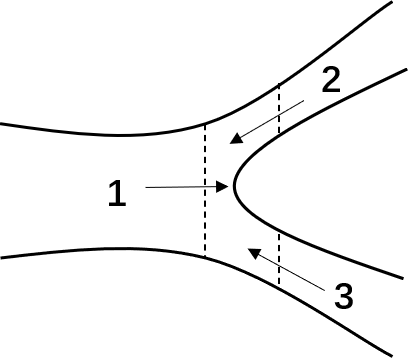}
  \caption{Illustration of the blood circulation at a vascular bifurcation.}
\end{figure}

There are three situations in the human body's vascular bifurcation structure: one in and multiple out, multi-segment confluence, and one in and one out (the parameters of the two connected segments are very different). This section will take the one-in-multiple-out situation as an example (as shown in \cref {fig:vessel_bif} shown), a detailed introduction to the model and processing of the vascular bifurcation \cite{Sherwin2003},\cite{Pedley1996}.

Based on the Riemann invariant method, we obtain the characteristic quantities of \cref{eq:1d_blood_set} and plug into the different tube law functions of arteries and veins $\phi (A(x,t),A_0(x),\beta (x))$, then we achieve the characteristic quantity of the artery
\begin{equation*}
  W_1 = v+4\sqrt{\frac{\beta}{2\rho}}A^{\frac{1}{4}}\qquad W_2 = v-4\sqrt{\frac{\beta}{2\rho}}A^{\frac{1}{4}},
\end{equation*}
and the characteristic quantity of the vein
\begin{equation*}
  W_1 = v+\int_{A_0}^{A}\frac{c(a)}{a}\mathrm{d}a\qquad W_2 = v-\int_{A_0}^{A}\frac{c(a)}{a}\mathrm{d}a.
\end{equation*}
As shown in \cref{fig:vessel_bif}, the forward wave of vessel 1 and the backward wave of vessels 2 and 3 together affect the blood flow distribution at the bifurcation. It is known that in each segment of the blood vessel, the characteristic quantity at the position $x$ remains unchanged. Combining the law of total pressure and the law of conservation of mass, we have the following equations of the blood flow characteristics at the bifurcation point:
\begin{equation}\label{eq:vessel_bif_char}
  \begin{aligned}
     & W_1^1 = v_1+4\sqrt{\frac{\beta_1}{2\rho}}A_1^{\frac{1}{4}} ,                       \\
     & W_2^2 = v_2-4\sqrt{\frac{\beta_2}{2\rho}}A_2^{\frac{1}{4}}  ,                      \\
     & W_3^2 = v_3-4\sqrt{\frac{\beta_3}{2\rho}}A_3^{\frac{1}{4}} ,                       \\
     & Q = v_1A_1 = v_2A_2 + v_3A_3  ,                                                    \\
     & P = \rho \frac{v_1^2}{2}+p_1 = \rho \frac{v_2^2}{2}+p_2 =\rho \frac{v_3^2}{2}+p_3,
  \end{aligned}
\end{equation}
the solutions of which are blood flows and velocities $(A_1,v_1), (A_2,v_2)$ and $(A_3,v_3)$ at the bifurcation point.
Specifically, using FEM, we can discrete the \cref{eq:vessel_bif_char} as
\begin{equation}\label{eq:alg_vessel_WQP}
  \begin{aligned}
     & W_1^1(A_{1,L_1}^{n+1},v_{1,L_1}^{n+1}) = W_1^1(A_{1,L_1-1}^{n},v_{1,L_1-1}^{n}) ,                     \\
     & W_2^2(A_{2,0}^{n+1},v_{2,0}^{n+1}) = W_2^2(A_{2,1}^{n},v_{2,1}^{n}) ,                      \\
     & A_{1,L_1}^{n+1}v_{1,L_1}^{n+1} = 2A_{2,0}^{n+1}v_{2,0}^{n+1} ,                       \\
     & p_1 + \frac{1}{2}\rho(v_{1,L_1}^{n+1})^2 = p_2 + \frac{1}{2}\rho(v_{2,0}^{n+1})^2;
  \end{aligned}
\end{equation}
so we can solve for $A_{1,L_1}^{n+1}$, $v_{1,L_1}^{n+1}$, $A_{2,0}^{n+1}$, $v_{2,0}^{n+1}$ using Newton iteration method.
Similarly, we can construct equations for the multi-stage confluence and one-in-one-out bifurcation form, respectively.

\subsection{Boundary Conditions}
When using the blood vessel model \cref{eq:1d_blood_set} to numerically simulate blood vessels, we need to set appropriate boundary conditions. At the start of the segment of blood vessel, we generally apply experimentally determined blood flow function $Q_{in}(t)$, blood pressure function $P_{in}(t)$, cross-sectional area $A(0,t)$ or flow velocity function $v(0,t)$. At the end of the blood vessel segment, in addition to the experimentally determined blood flow function $Q_{out}(t)$, blood pressure function $P_{out}(t)$, cross-sectional area $A(l,t)$ or flow velocity function $v(l,t)$, there are the following three cases.

\subsubsection{Absorbing boundary.}

When the characteristic quantity of the blood vessel terminal is a constant independent of time, we call this boundary condition an absorbing boundary. According to \cref{eq:vessel_bif_char} under the absorbing boundary condition, the cross-sectional area $A(l,t)$ of the end of the blood vessel and the flow velocity function $v(l,t)$ are constant values independent of time. Therefore, at the end of the blood vessel, the cross-sectional area  $A(l,t)$ and the flow velocity function $v(l,t)$ can be set based on the reference cross-sectional area $A(0,t)$ and the reference flow velocity function $v(0,t)$, respectively.

\subsubsection{Coupling a terminal resistance.}\label{ch:terminal_couple}

\begin{figure}[htbp]
  \centering
\label{fig:terminal_couple}
\includegraphics{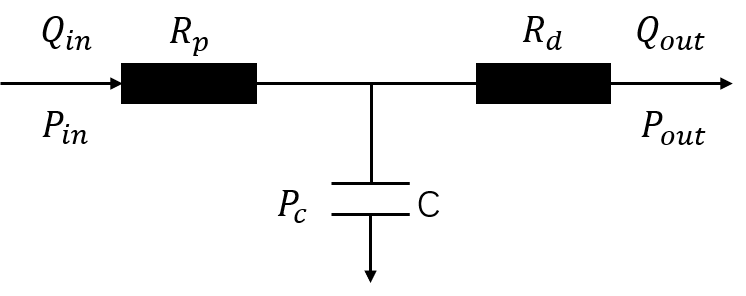}
  \caption{Illustration of Coupling a terminal resistance.}
\end{figure}

In the human body, the ends of arteries are often connected to capillaries, and blood flows through the capillaries to exchange materials with surrounding tissues and then flows into the veins. The resistance of peripheral blood vessels (capillaries) has an important influence on the blood pressure of the human body. Therefore, we couple a terminal resistance at the end of the arteries and blood vessels, as shown in Figure \cref{fig:terminal_couple}, to simulate the peripheral resistance of the human blood vessels.

Suppose the resistance $R$ represents the total peripheral resistance of the surrounding blood vessels, $Q_{in}$ and $Q_{out}$ are the blood flow from the arterial system and out of the capillaries, respectively, $P_{in}$ and $P_{ out}$ refer to the blood pressure at the end of the blood vessel and the right end of the circuit. The time dynamics of the cross section area and velocity $(A_*^{t+\Delta t},v_*^{t+\Delta t})$ satisfy
\begin{equation}\label{eq:bound_FEM_A}
A_*^{t+\Delta t}v_*^{t+\Delta t} = \frac{P(A_*^{t+\Delta t})-P_{out}}{R}.
\end{equation}

Due to small time and space steps in FEM, the Riemann features at adjacent positions are approximately equal, in other words,
\begin{align}\label{eq:bound_FEM_1}
 W_1{(A_*^{t+\Delta t},v_*^{t+\Delta t})} = W_1(A_l^{t},v_l^{t}),
\end{align}
\begin{align}\label{eq:bound_FEM_2}
W_2{(A_*^{t+\Delta t},v_*^{t+\Delta t})} = W_2(A_r^{t},v_r^{t}).
\end{align}

The cross-section area on the left and right to one vessel segment should be the same, that is, $A_l^{t} = A_r^{t}$. According to \cref{eq:1d_blood_relation}, \cref{eq:1D_blood_velocity}, \cref{eq:vessel_bif_char} and \cref{eq:bound_FEM_1}, \cref{eq:bound_FEM_2}, we have
\begin{equation}\label{eq:bound_FEM}
RA_*^{t+\Delta t}(v_l^{t}+4c(A_l^{t})-4c(A_*^{t+\Delta t})) = \frac{\beta}{A_0}(\sqrt{A_*^{t+\Delta t}}-\sqrt{A_0})+P_0-P_{out}
\end{equation}
With the known $(A_l^{t_0},v_l^{t_0})$ at the starting time $t=0$, we can use Newton's iteration method to solve for $A_*^{t_0+\Delta t}$ with \cref{eq:bound_FEM}, and $v_*^{t_0+\Delta t}$ by substituting $A_*^{t_0+\Delta t}$ into \cref{eq:bound_FEM_A}. Then through \cref{eq:bound_FEM} and the assumption of $A_l^{t} = A_r^{t}$ we can solve $(A_r^{t_0},v_r^{t_0})$.

Using the FEM of the 1-D blood vessel model in \cref{ch:FEM_1D_vessel}, we can find $A(l-\Delta x,t_0)$ and $v(l-\Delta x,t_0)$ with iteration. Next, based on the calculated $(A_r^{t_0},v_r^{t_0})$, we can use the aforementioned FEM to calculate $(A_l^{t_0+\Delta t},v_l^{t_0+\Delta t}) $.
Similarly, starting from the moment $t = t_0 + \Delta t$, we can obtain the value of $(A_l^{t},v_l^{t})$ at each moment using iteration, which is the boundary condition at the end of the blood vessel .

\subsubsection{Coupling an RCR circuit}\label{ch:RCR_couple}

\begin{figure}[htbp]
  \centering
\label{fig:RCR_couple}
\includegraphics{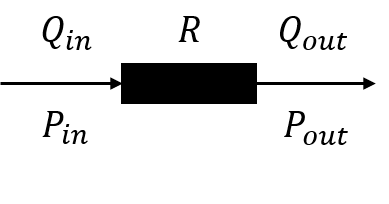}
  \caption{Illustration of Coupling a RCR circuit.}
\end{figure}

The human body's blood vessels can have a certain degree of compliance, that is, the blood vessels can increase their volume without breaking under the pressure. To reflect the peripheral vascular resistance and the compliance characteristics of the blood vessels simultaneously, when the blood vessel terminal is connected to the capillary, we consider coupling an RCR circut to the end of the blood vessel. Specifically, the RCR model is a lumped parameter circuit model, including a remote resistance $R_p$, a near end resistance $R_d$, and a capacitor $C$, as illustrated in \cref{fig:RCR_couple}.

Suppose the total peripheral resistance of the blood vessel is $R$, i.e. $R_p + R_d = R$, $Q_{in}$ and $Q_{out}$ represent the blood flow from the arterial system and out of the capillaries; $P_{in} $, $P_{out}$ and $P_C$ represent blood pressure at the blood vessel terminal, circuit terminal and capacitor $C$, respectively. We show how to solve the vascular terminal cross-sectional area $A(l,t)$ and blood flow velocity $v(l,t)$ at time $t$ in an RCR circuit. First, $P_C$ satisfies
\begin{equation*}
C\frac{P_C}{\mathrm{d}t} = A_*^{t}v_*^{t} - \frac{P_C-P_{out}}{R_d},
\end{equation*}
and the first differentiation of $P_C$ in terms of time $t$ becomes
\begin{equation}\label{eq:RCR_PC}
P_C^n = P_C^{n-1} + \frac{\Delta t(A_l^{(n-1)\Delta t}v_l^{(n-1)\Delta t}-(P_C^{n-1}-P_{out})/R_d)}{C}.
\end{equation}
$P_C^0 = 0$ is the initial blood pressure of the venous system for $P_{out}$. Using \cref{eq:RCR_PC}, we calculate $P_C$. Similarly, we can derive
\begin{equation}\label{eq:alg_RA_solve_Av}
\begin{aligned}
R_pA_*^{t+\Delta t}(v_l^{t}+4c(A_l^{t})-4c(A_*^{t+\Delta t})) = \frac{\beta}{A_0}(\sqrt{A_*^{t+\Delta t}}-\sqrt{A_0})+P_0-P_{C},\\
A_*^{t+\Delta t}v_*^{t+\Delta t} = \frac{P(A_*^{t+\Delta t})-P_{C}}{R_p},
\end{aligned}
\end{equation}

which can result in values of other parameters using the idea in \cref{ch:terminal_couple}.

\subsection{Modeling a liver}\label{ch:liver_model}

\begin{figure}[htbp]
  \centering
\label{fig:liver}
\includegraphics[scale = 0.25]{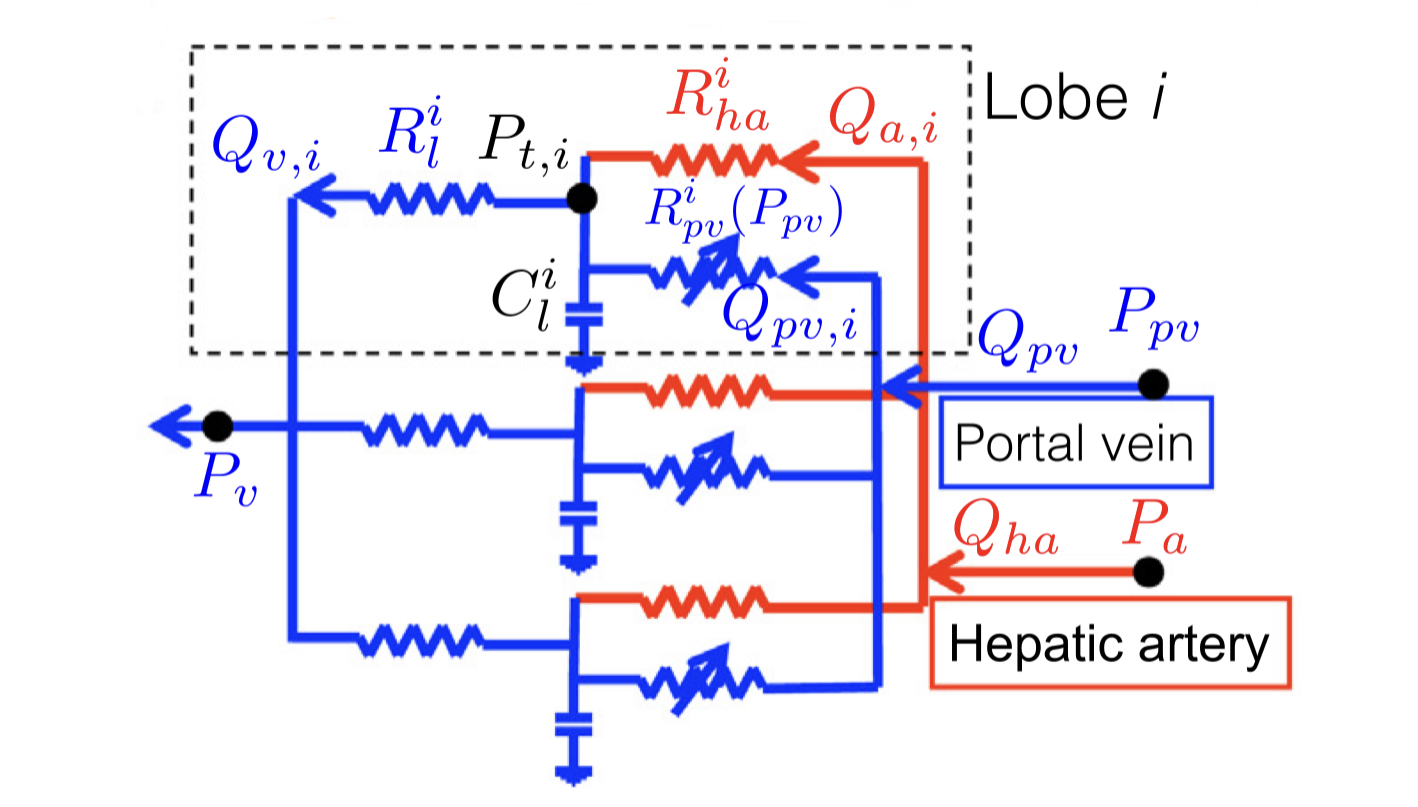}
  \caption{Modeling a human liver.}
\end{figure}

 The human liver is made of three independent liver lobes covered with a blood vessel network. The main components of the hepatic artery (HA) and portal vein (PV) provide the liver with blood sources, while the hepatic vein (HV) leads the blood flow in the liver out of the liver and into the inferior vena cava. This two-in-and-one-out vascular structure reflects the specificity of the blood supply of the liver. The blood supply system of the human liver is similar to that of the pig liver, so our liver model bases on the research on pig liver \cite{Audebert2017} and \cite{Audebert2017partial}. Here, we simulate the three main liver lobes in the liver using three parallel circuits, each circuit contains the flow input from the hepatic artery (HA) and the portal vein (PV), and the flow output from the hepatic vein (HV), as illustrated in \cref{fig:liver}.

Assuming that the resistance value in the circuit model corresponding to each part of the liver lobe is proportional to the mass of that part of the liver lobe, and the capacitance is inversely proportional to the mass of this part of the liver lobe, we have  for each part of the liver lobe: $i = 1, 2, 3$,
\begin{align}
    & P_a-P_{t,i}=R_{ha}\frac{M}{M_i}Q_{a,i}, \label{eq:liver_system_1}       \\
    & P_{pv}-P_{t,i}=R_{pv}\frac{M}{M_i}Q_{pv,i},  \label{eq:liver_system_2}  \\
    & P_{t,i}-P_v=R_l\frac{M}{M_i}Q_{v,i},    \label{eq:liver_system_3}       \\
    & C_lM_i\frac{dP_t}{dt}=Q_{a,i}+Q_{pv,i}-Q_{v,i}, \label{eq:liver_system_4}
\end{align}

Based on \cref{eq:liver_system_1}, \cref{eq:liver_system_2} and \cref{eq:liver_system_3}, we can write $Q_{a,i}$,$Q_{pv,i}$, $Q_{v,i}$ as an expression of $P_t$, and then substitute them into \cref{eq:liver_system_4}, which then can be transformed into an ordinary differential equation about $P_t$. We can then derive the analytical solution of $P_t$ as
\begin{equation}\label{eq:liver_Pt}
P_t = (P_0 + \frac{C_1}{C_2})e^{C_2t}+\frac{C_1}{C_2},
\end{equation}
where $P_0$ is the initial blood pressure of the liver,
\begin{align*}
& C_1 = \frac{1}{C_lM}(\frac{P_a}{R_{ha}}+\frac{P_{pv}}{R_{pv}}-\frac{P_v}{R_l})
& C_2 = -\frac{1}{C_lM}(\frac{1}{R_{ha}}+\frac{1}{R_{pv}}-\frac{1}{R_l})
\end{align*}

To estimate parameters in the liver model, we first re-write \cref{eq:liver_system_1}, \cref{eq:liver_system_2} and \cref{eq:liver_system_3} as
 \begin{equation*}
      \left\{
\begin{aligned}
    & Q_{a,i}=\frac{M_i\left(P_a-P_t\right)}{MR_{ha}},    \\
    & Q_{pv,i}=\frac{M_i\left(P_{pv}-P_t\right)}{MR_{pv}}, \\
    & Q_{v,i}=\frac{M_i\left(P_t-P_v\right)}{MR_l},
    \end{aligned}
      \right. 
\end{equation*}
The sum of which across $i =1,2,3$ is
 \begin{equation}\label{eq:liver_3Qs}
      \left\{
     \begin{aligned}
         & Q_a=\frac{P_a-P_t}{R_{ha}}       \\
         & Q_{pv}=\frac{P_{pv}-P_t}{R_{pv}} \\
         & Q_v=\frac{P_t-P_v}{R_l}
     \end{aligned}
      \right.
\end{equation}
Hence, we can calibrate the resistances of HA, PV and HV through
\begin{equation}\label{eq:liver_resistance}
      \left\{
      \begin{aligned}
         & R_{ha}=\frac{P_a-P_t}{Q_a}       \\
         & R_{pv}=\frac{P_{pv}-P_t}{Q_{pv}} \\
         & R_l=\frac{P_t-P_v}{Q_v},
      \end{aligned}
      \right.
\end{equation}
the sum of the last equation above across $i = 1, 2, 3$ leads to
\begin{equation}\label{eq:liver_C_ODE}
      C_lM\frac{dP_t}{dt}=Q_a+Q_{pv}-Q_v
\end{equation}

Because that the blood inflow and outflow an internal organ are believed to be equal \cite{Markou2004}, that is, $Q_v=Q_a+Q_{pv}$. \cref{eq:liver_C_ODE}  indicates that $\frac{dP_t}{dt }$ is approximately equal to 0, so we treat hepatic lobe blood pressure $P_t$ as a constant approximately.
Because the blood pressure and blood flow of a human liver cannot be directly obtained during normal physiological activities using invasive measuring, we consider using adjacent blood vessels in our numerical study.
As a branch of the abdominal aorta (ABO), we can examine the HA blood , namely the blood pressure $P_a$ and blood flow $Q_a$ at the entry of the branch, via ABO and liver branches using a 1D blood vessel model.
Because the veins are not affected by pulse waves and the blood flow resistance of the venous blood vessels is small, the blood pressure of the veins in a human body changes little and the blood pressure remains at a low level.
The blood pressure of the human right atrium and large veins in the thoracic cavity is usually called central venous pressure (CVP), the value of which is relatively low and ranging $0.4-1.2 kPa$ \cite{Perthame2002}. In our numerical study, we use the estimated CVP value as hepatic PV blood pressure $P_{pv}$.
The HV is a branch of the inferior vena cava (IVC), hence we use the estimated IVC blood pressure as the HV blood pressure $P_v$\cite{Audebert2017partial}.
The liver blood pressure $P_t$ is between $P_v$ and $P_{pv}$. Generally the blood pressure drop between the PV and the liver accounts for 80\% of the blood pressure drop between the PV and the HV \cite{Audebert2017partial}.
Moreover, the blood supply ratio of the human HA to the PV is about 3:1 \cite{Sherwin2003a}. Therefore, after estimating the HA blood flow  $Q_a$, we can calculate the PV and HV blood flow $Q_{pv}$ and $Q_v$ accordingly.

\section{Numerical Study}
\label{sec:alg}

\subsection{Data Description}

Our experimental data bases on a sample of 33 healthy humans under normal physiological conditions. The sample age ranges from 20 to 50 years old. The doctor collected each individual's measurements regarding the vascular and blood flow characteristics below. \begin{itemize}
    \item \textbf{PS:} Peak systolic flow velocity, that is, the peak blood flow velocity during systole;
    \item \textbf{ED:} The lowest flow velocity during diastole, that is, the minimum blood flow velocity during diastole. Note that there may be back-flow of blood during diastole
    \item \textbf{VTI:} Blood flow velocity time integral, that is, the area enclosed by the blood flow velocity spectrum and the baseline or time axis in a cardiac cycle, which is equivalent to the height of an ideal liquid column;
    \item \textbf{Heart rate:} The number of people's heartbeats per minute;
    \item \textbf{Diameter of blood vessel;}
    \item  \textbf{Blood flow:}  The volume of blood flow through the cross section of the blood vessel per unit time.
\end{itemize}

\begin{figure}[htbp]
  \centering
\label{fig:flow_distribution}
\includegraphics[scale = 0.4]{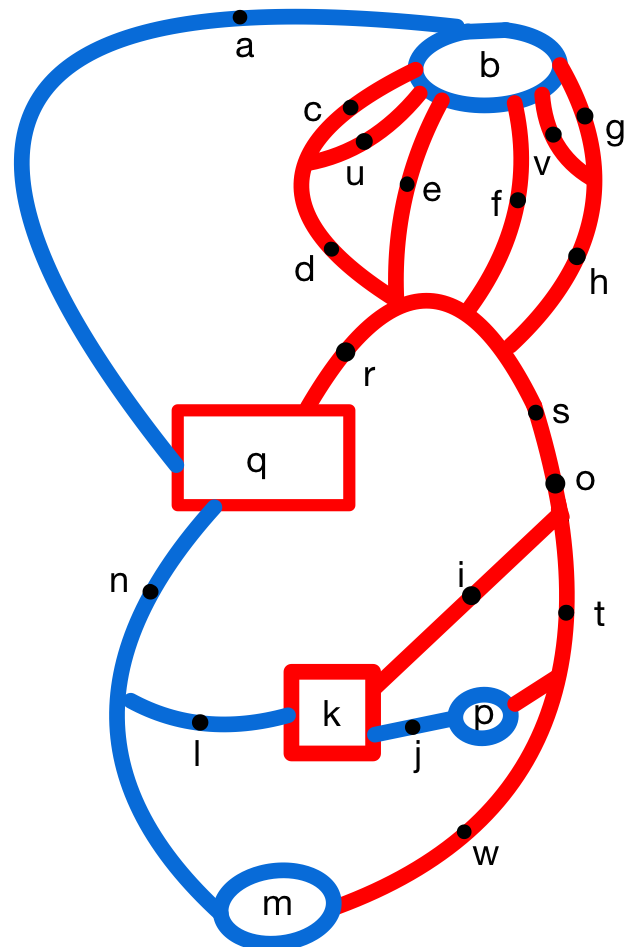}
  \caption{Illustrative model on the human blood circulation.}
\end{figure}

We show the schematic diagram of the human blood circulation used in our model in \cref{fig:flow_distribution}, and the names of the observed points on blood vessels and organs in \cref{fig:flow_distribution} are shown in  \cref{tab:names_circulation}.

The measurements are taken at abdominal aorta A (point o), right common carotid artery (point e), and left common carotid artery (f Point), right brachial artery (point d), left brachial artery (point h), right radial artery (point c), left radial artery (point g).

\begin{table}[tbhp]
{\footnotesize
  \caption{Names of organs and observed points in human blood circulation}\label{tab:names_circulation}
\begin{center}
  \begin{tabular}{cc|cc} \hline
   Position     & Name            & Position       & Name     \\ \hline
    a      & Superior vena cava     & m  & surrounding capillary network         \\
    b & capillary network    & n  & inferior vena cava    \\
    $c^*$  & right radial artery & $o^*$ & abdominal aorta A \\
    $d^*$ & right brachial artery  & p  & digestive system capillary network    \\
    $e^*$ & right common carotid artery         & q   & heart     \\
    $f^*$  & left common carotid artery     &  r       &    ascending aorta        \\
    $g^*$  & left radial artery     &  s       &     descending aorta         \\
    $h^*$  & left brachial artery      &  t       &    abdominal aorta B         \\
    i  & hepatic artery     &  u       &    right ulnar artery      \\
    j  & portal vein      &  v       &    left ulnar artery         \\
    k  & liver    &  w       &    abdominal aorta C      \\
    l  & hepatic vein      &        &             \\
    \hline
  \end{tabular}
\end{center}
}
\end{table}

\subsection{Algorithm to estimate parameters in blood vessels}

\begin{table}[tbhp]
{\footnotesize
  \caption{Table of parameters}\label{tab:paremeter}
\begin{center}
  \begin{tabular}{lll} \hline
  \multicolumn{3}{l}{Regarding blood vessels} \\ \hline
   Parameter      & Name        & Unit    \\ \hline
     $\beta$       & Blood vessel elasticity      & $dyn/cm^3$    \\
    $\rho $       & Blood density      & $g/cm^3$      \\
    $\mu$         & Blood viscosity    & $dyn/cm^2$    \\
    $K_f$         & Friction coefficient  & $cm^2/s$             \\
    $R_d$         & Remote resistance      & $dyn\cdot s/cm^5$ \\
    $R_p$         & Near-end resistance  & $dyn\cdot s/cm^5$ \\
    $C$           & Capacitance          & $cm^5/dyn$        \\
\hline
  \end{tabular}
\end{center}
}
\end{table}
Please note that the blood density, blood viscosity and friction coefficient have relatively small differences among individuals, hence we can set them constants that do not change with individuals. Vascular elasticity parameters are simultaneously affected by the individual's age and physiological condition, and thus can vary with individuals. So we need to estimate it separately for each individual.

In the  diagram of the systemic circulation \cref{fig:flow_distribution}, the radial artery and the ulnar artery are two branches of the brachial artery. Because the diameter of the radial artery and the ulnar artery are relatively small, for simplicity, here we consider the radial artery and the ulnar artery as two identical branches of the brachial artery. Specifically, we use the experimental data of the $d-c$ and the $h-g$ segments to get the approximate value of the blood vessel elasticity parameter in \cref{alg:blood}.

The first part of the algorithm aims to estimate parameters of different individuals' blood vessels $\beta_1$ using observed left radial artery blood flow $Q$, cross-section areas and velocities at the observed points $h$ and $g$: $v_1(0,t)$, $A_1(0,t)$, $v_2(l_2,t)$, $A_2(l_2,t)$, $v_1(x,0)$, $A_1(x,0)$, $v_2(x,0)$ and $A_2(x,0)$. The reference cross-section area at the two observed points $h$ and $g$ are $A_{01}$, $A_{02}$.
Then we estimate each individual's blood vessel elasticity coefficient $\beta_2$ using the data of the right brachial artery-right radial artery ($d-c$ in \cref{fig:flow_distribution}). To reduce the error, we take the mean value of $\beta_1$ and $\beta_2$ as the estimated value of the blood vessel elasticity coefficient $beta$ for the sample of individuals.

The second part of the algorithm focusing on estimating the resistance in the RCR circuit model. Because the capacitance parameters in the RCR model have almost no individual specificity, we can set each individual's capacitance parameter $C$ to a uniform constant.
We can calibrate the total resistance $R_{tot}$ in an RCR circuit using the formula $R_{tot} = \frac {\Delta P}{Q}$, where $\Delta P$ represents the difference between the blood pressure at the end of the arterial vessel and the initial end of the venous vessel. Without loss of generality, we can assume that the venous blood pressure is zero, and $Q$ represents the blood flow into the vein through the arterial vessel. Moreover, the remote-end resistance $R_d$ and the near-end resistance $R_p$ account for $90\%$ and $10\%$ \cite{Audebert2017} of the total resistance, respectively. Assuming that the same RCR model is coupled to the end of the arteries connecting capillaries in the same human body, we can find the resistance value in the RCR circuit by calculating the end of the arterial blood pressure from the blood pressure of the end of the left and right radial arteries (observed points $c$ and $g$).

The space step is $\Delta x$ and time step is $\Delta t$.
Suppose space points $i = 0, 1, 2, \cdots, l/\Delta x$ and time points $n = 0, 1, 2, \cdots, T/\Delta t$.

\begin{algorithm}
\caption{Estimate parameters of blood vessels }
\label{alg:blood}
\begin{algorithmic}
\FOR{Blood vessel segments $h-g$ and $c-d$}
\STATE{Initiate $A_{q,i}^n$, $v_{q,i}^n$, $q = 1, 2$ at the pair of observed points ($h$ and $g$, or $c$ and $d$.}
\WHILE{$1\cdot 10^5 \leq \beta \leq 1\cdot 10^6 dyn/cm^3$}
\STATE{Update blood parameters using \cref{alg:update_AXv} when $q = 1$.}
\ENDWHILE
\STATE{Use Newton's iteration method to solve \cref{eq:alg_vessel_WQP}. }
\STATE{Update blood parameters using \cref{alg:update_AXv} when $q = 2$.}
\STATE{Calculate the blood flow at the midpoint of the vessel segment $Q(\frac{l_2}{2})= \frac{1}{T+1}\sum_{n=0}^{n=T} A_{2,\lfloor \frac{L_2}{2}\rfloor}^nv_{2,\lfloor \frac{L_2}{2}\rfloor}^n$, $DIFF_{\beta} = |Q-Q(\frac{l_2}{2})|$.}
\STATE{Find $\beta = arg\mathop{min}\limits_{\beta} DIFF_{\beta}$}
\ENDFOR
\RETURN the average value of $\beta$.
\STATE{Initiate length of the vessel $l=10cm$, reference blood pressure $P_0=10^5 dyn/cm^2$, density $\rho=1.05 g/cm^3$ viscosity $\mu=0.04 dyn/cm^3$, friction $K_f=22\pi \frac{\mu}{\rho}$,  each individual's $\beta$; $A_{1,i}^n$,  $v_{1,i}^n$ at the observed points ($h$ and $c$).}
\STATE{Update blood parameters using \cref{alg:update_AXv} when $q = 1$.}
\FOR{Each $n$}
\STATE{Calculate blood pressure $P^n$ according to \cref{eq:alg_P_update}. }
\ENDFOR
\STATE{$P = \frac{1}{T+1}\sum_{n=0}^{T}P^n$, $R_{tot} = \frac{P}{Q}$.}
\RETURN $R_p = 0.1R_{tot}$, $R_d = 0.9R_{tot}$.
\end{algorithmic}
\end{algorithm}

\subsection{Algorithm on liver parameters}

To estimate the parameters in the liver model, we first obtain the blood flow $Q_a$ and blood pressure $P_a$ at the initial site of the HA, through a 1D blood vessel model on the abdominal aorta (ABO)-HA bifurcation vessel. In \cref{fig:flow_distribution}, ABO bifurcates below $o$ point---a part of the blood flows into the liver and a part of the blood flows into the abdominal aorta B (ABO-B) segment. Hence, we model this ABO-HA bifurcation with the one-in-and-two-out model in \cref{ch:vessel_bif_model}. Specifically, the end of abdominal aorta A (position $o$) bifurcates into the beginning of the HA (position $i$) and abdominal aorta B (position $t$). We assume the cross-sectional area of abdominal aorta B (position $t$) is one-half of that of abdominal aorta A (position $o$) \cite{barrett2010ganong}, and the blood flow velocity is zero at the absorption boundary. The end of the HA (position $i$) connects to the capillary, which can be regarded as a coupled RCR model in \cref{ch:RCR_couple}.

The first part of \cref{alg:liver} solves the blood flow velocity and cross-sectional area at the end of the HA (blood segment $i$), using the measurement at the abdominal aorta A vessel segment (position $o$).
Though we cannot directly measure the cross-sectional area of an HA, medical literature reveals that the diameter of the HA is generally $2-5mm$ \cite{rubin1998measurement}, which is similar to the diameter of the brachial artery. Thus, we use the measured cross-sectional area of the brachial artery (position $d$ or $h$???) as a proxy of the HA cross-sectional area. Based on the measured experimental data, namely the PS, ED and the blood flow velocity of the abdominal aorta A (position $o$), we set up the blood flow velocity change in a cardiac cycle of length $Period$ for each individual, according to he proportion of the various phases of the human cardiac cycle in the cardiac cycle and the law of blood flow changes \cite{barrett2010ganong}. Typically the first $0.15\,Period$ in a cardiac cycle, the blood velocity is zero, and then it increases to the peak during $0.15-0.25\,Period$, followed by decreasing to the lowest level between $0.25$ and $0.42 \,Period$. When the blood velocity raises back to zero from $0.42$ to $0.45\,Period$ and remains zero until the end of the cycle.
In \cref{alg:liver}, we take the time interval $\Delta t = 0.0001s$ and the space interval $\Delta x = 0.1cm$.
Let the vessel segments $l_1$, $l_2$ and $l_3$ be the abdominal aorta A vessel marked with $o$, the HA labeled with $i$ and abdominal aorta B at $t$ of \cref{fig:flow_distribution}, respectively.
In addition to the blood parameters, namely vascular elasticity coefficient $\beta$, remote resistance $R_d$ and near-end resistance $R_p$ estimated in \cref{alg:blood}, \cref{tab:alg_input_liver} presents input values from literature required for \cref{alg:liver}.

\begin{figure}[htbp]
  \centering
\label{fig:blood_vol_cycle_ABOA}
\includegraphics[scale = 0.4]{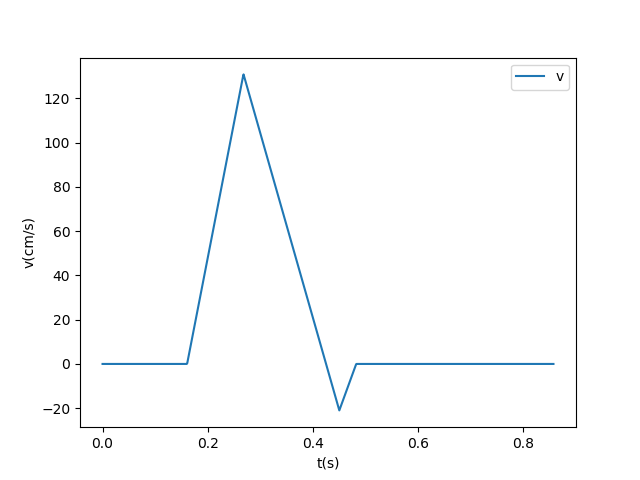}
  \caption{The blood velocity changes at abdominal aorta A (position $o$ in a cardiac cycle.}
\end{figure}

\begin{table}[tbhp]
{\footnotesize
  \caption{Table of inputs for \cref{alg:liver}}\label{tab:alg_input_liver}
\begin{center}
  \begin{tabular}{llll} \hline
  \hline
   Parameter      & Abdominal         & HA & Abdominal    \\
    &     aorta A &   & aorta B \\\hline
    Blood vessel length $l_i$ \cite{rubin1998measurement}& $5cm$ & $5cm$ & $5cm$ \\
    Reference blood pressure $P_0$ \cite{Audebert2017} & $10^5 dyn/cm^2$ & $10^5 dyn/cm^2$ & $10^5 dyn/cm^2$ \\
    Initial blood velocity  $v_(i)(x,0)$ \cite{Audebert2017} & $0cm/s$ & $0cm/s$ & $0cm/s$ \\
    Blood density$\rho$ \cite{martini2006human}& $1.06 g/cm^3$ & $1.06 g/cm^3$ & $1.06 g/cm^3$ \\
    Blood viscosity coefficient  $\mu$ \cite{Audebert2017} & $0.04 dyn/cm^3$ & $0.04 dyn/cm^3$ & $0.04 dyn/cm^3$ \\
    Friction coefficient $K_f$ \cite{Audebert2017} & $22\pi \frac{\mu}{\rho}$ & $22\pi \frac{\mu}{\rho}$ & $22\pi \frac{\mu}{\rho}$ \\
        Capacitance  $C$ \cite{Audebert2017} &-& $4\cdot 10^{-5}cm^5/dyn$ &-\\
\hline
  \end{tabular}
\end{center}
}
\end{table}

The second part of \cref{alg:liver} aims to estimate the liver model parameters for each individual in the sample based on \cref{ch:liver_model}, namely the resistance parameters of the HA, PV and HV,  $R_{ha}$, $R_{pv}$ and $R_l$. First, according to $Q_v=Q_a+Q_{pv}$ and $Q_{pv}=3Q_a$, we can obtain the blood flow in the PV and HV  ($Q_{pv}(ml/s)$ and $Q_v(ml/s)$).
Human physiology finds that the blood pressure value can drop to $2-2.7kPa$ when blood flows through arteries and capillaries to venules in a human's blood circulation, and the blood pressure is about $0.4 -0.5kPa$ when flowing to the inferior vena cava \cite{Markou2004}.
In \cref{alg:liver}, we use the estimated value of central venous pressure  as the average PV blood pressure $P_{pv}$, and the estimated value of the inferior vena cava blood pressure as the average PV blood pressure $P_v$. Finally, due to $P_{pv}-P_t = 80\% (P_{pv}-P_v)$ \cite{Audebert2017partial}, we can obtain the estimated value of the average liver blood pressure $P_t$.

\begin{algorithm}
\caption{Estimate parameters in the liver model }
\label{alg:liver}
\begin{algorithmic}
\STATE{Initiate known parameters in \cref{tab:alg_input_liver}; initial boundary value conditions $v_1(0,t)$, $A_1(0,t)$, $v_3(l_3,t)$, $A_3(l_3,t)$, $v_i(x,0) $,$A_i(x,0)$($i=1,2,3$); space step size $\Delta x$; time step size $\Delta t$; reference cross-sectional area $A_{0i}$, $L_i = \lfloor \frac{l_i}{\Delta x}\rfloor$, $i=1,2,3$;}
\FOR{$i = 1, 2, 3$}
 \WHILE{$j = 0 \quad to \quad L_i$ }
\STATE{Calculate $A_{i,j}^0$, $v_{i,j}^0$ based on initial inputs}
\ENDWHILE
\ENDFOR
\WHILE{$n = 0 \quad to \quad T$ }
 \STATE{Set up $A_{1,0}^n$, $v_{1,0}^n$, $A_{3,L_3}^n$, $v_{3,L3}^n$ using initial inputs.}
\ENDWHILE
\WHILE{$n = 0 \quad to \quad T$ }
  \WHILE{$i = 1 \quad to \quad L_1 - 1$ }
    \STATE{Update blood parameters using \cref{alg:update_AXv}.}
  \ENDWHILE
  \STATE{Use Newton's iteration method to solve \cref{eq:alg_vessel_WQP}. }
  \FOR{i = 2, 3}
   \WHILE{$j = 1 \quad to \quad L_i - 1$ }
     \STATE{Update blood parameters using \cref{alg:update_AXv}.}
   \ENDWHILE
  \ENDFOR
  \STATE{Update $P^{n+1}_C using \cref{eq:RCR_PC}$;}
  \STATE{Solve for $A_*$, $v_*$ in \cref{eq:alg_RA_solve_Av};}
  \STATE{Let $A_r = A_{2,L_2}^{n}$, solve for $v_r$ using the equation $W_2(A _*,v_*) = W_2(A_r,v_r)$ derived from \cref{eq:bound_FEM};}
   \STATE{$A_{2,L_2+1}^n = A_r$, $v_{2,L_2+1}^n = v_r$, update $A_{2,L_2}^{n+1}$, $v_{2,L_2}^{n+1}$ using \cref{alg:update_AXv}}
 \ENDWHILE
\WHILE{$n = 0 \quad to \quad T$}
  \STATE{Update $P^n$ according to \cref{eq:alg_P_update}, $Q^n = A_{2,0}^nv_{2,0}^n$; }
\ENDWHILE
\RETURN Boundary values at the end of HA: $P = \frac{1}{T+1}\sum_{n=0}^{T}P^n$, $Q = \frac{1}{T+1}\sum_{n=0}^{T}Q^n$.
\STATE{Calculate the blood flow in the PV and HV  ($Q_{pv}(ml/s)$ and $Q_v(ml/s)$ using $Q_v=Q_a+Q_{pv}$ and $Q_{pv}=3Q_a$;}
\STATE{Estimate $P_t$ using $P_{pv}-P_t = 80\% (P_{pv}-P_v)$;}
\STATE{Solve for $R_{ha}$, $R_{pv}$ and $R_l$ using \cref{eq:liver_resistance}}
\RETURN $R_{ha}$, $R_{pv}$ and $R_l$.
\end{algorithmic}
\end{algorithm}

\subsection{Verify numerical results}
\label{ch:verify_results}

We verify the results from \cref{alg:blood} and \cref{alg:liver} by comparing the estimated values with the observed measurements.
Because it is impossible to directly measure the relevant parameters of the human internal organs through invasive clinical methods, we estimate each individual's blood inflow into the heart through the vein, and compare it with the blood outflow of the heart measured in the clinical experiment. We then examine the mean relative error $MRE$ between the estimated inflow and measured outflow across the individuals in the sample.
Specifically, the heart pumping blood volume at position  $r$ is actually the blood flow of the ascending aorta at position $a$---the superior vena cava, which is the sum of and the blood volume at observed positions $d$, $e$, $f$, $h$ and  $o$ \cite{Fitzgerald2018}.
Moreover, the blood inflow into the heart through the vein is the sum of the blood flow of the superior vena cava at position $a$ and the inferior vena cava  at position $n$ in \cref{fig:flow_distribution}.
We describe the method to estimate the blood flow of the inferior vena cava at position $n$ in \cref{alg:inferiorVC}, where we substitute the liver parameters resulted from \cref{alg:liver} into the liver circuit model (\cref{ch:liver_model}), and calculate the blood flow of the HA (position $i$), of the PV at position $j$ and of the HV at position $l$.
Furthermore, the blood flow of abdomen aorta C (position $c$) is equal to the blood flow of abdominal aorta A (position $o$) minus the blood flow through the HA (position $i$) and the PV (position $j$). Therefore, the blood flow of the inferior vena cava (position $n$) is the sum of the blood flow through the HV (position $l$) and the abdominal aorta C ((position $w$).

\begin{algorithm}
\caption{Estimate blood flow in the inferior vena cava  }
\label{alg:inferiorVC}
\begin{algorithmic}
\STATE{Initiate liver resistance parameters estimated in \cref{alg:liver}, capacitance parameter $C_l=1.5cm^5/dyn$ \cite{Audebert2017}, cardiac cycle duration $Period$, hepatic artery blood pressure in one cardiac cycle $P_a(t)$; portal blood pressure $P_{pv}(t)$, liver Venous blood pressure $P_v(t)$, blood flow at position $o$ $Q_o$, time step $\Delta t$; $T = \lfloor \frac{Period}{\Delta t}\rfloor$;}
\FOR{$n = 0, 1, 2, \cdots, T$}
\STATE{Update liver blood pressure $P_{n\delta t}$ based on \cref{eq:liver_Pt};}
 \STATE{Calculate blood flow in the HA, PV and PA according to \cref{eq:liver_3Qs};}
\ENDFOR
\STATE{$Q_a = \frac{1}{T+1}\sum_{n=0}^{T}Q_a^n$, $Q_{pv} = \frac{1}{T+1}\sum_{n=0}^{T}Q_{pv}^n$, $Q_{v} = \frac{1}{T+1}\sum_{n=0}^{T}Q_{v}^n$;}
\RETURN Blood flow at position $n$: $Q_n = Q_o - Q_a-Q_{pv}+Q_{v}$.
\end{algorithmic}
\end{algorithm}

We obtain blood outflow from the heart from \cref{alg:inferiorVC}. Comparing with the observed blood inflow into the heart across the individuals in the sample, the mean relative error (MRE) is $1.3\%$.
Under normal physiological conditions, spontaneous breathing may lead to changes in blood pressure in the thoracic cavity. Therefore, the arterial pulse pressure decreases during inhalation and increases during exhalation, which leads to a change of heart stroke bleeding in each cardiac cycle. The differences in the amount of heart stroke bleeding are generally quantified with the stroke volume variation (SVV) \cite{34}, in other words, $SVV = \frac{SV_{max}-SV_{min}}{SV_{mean}}$, where $SV_{max}, SV_{min}$ and $SV_{mean}$ are the maximum, minimum and average stroke volume of the heart stroke bleed during a period of time. In  general,  the range of $SVV$ should be within $10\%$ as . From the perspective of stroke volume variation, $MRE <10\%$ indicates  that our estimated results of liver parameters are reasonable.

\section{Conclusions}
\label{sec:conclusions}

In this study, we use the blood vessel model as an important connection model to construct blood vessels that circulate between the heart and the liver in human blood circulation system.
Based on the study of blood flow in the blood vessel network, heart modeling, and liver modeling, this paper combines the three to form a circulatory blood propagation model of the heart, Vessels and liver in human bodies.
Our work features using non-invasive clinical readings to estimate individual-specific blood parameters.
required for the three-part model.
We contribute to clinical literature with
hemodynamic equations and finite element methods, which can effectively simulate the propagation of blood in arteries and veins. By analyzing the existence and stability conditions, we select appropriate time and space intervals to make the model solvable.
Our approach gives examples of several commonly used boundary conditions, explains the application scenarios of various boundary conditions, and derives the calculation process of boundary conditions in details. When the blood vessel or the connection state of the blood vessel changes, this approach requires only a slight modification of the equation structure and related parameters, without changing the numerical solution method of the model. That is, the model has strong adaptability in other scenarios.
To estimate model parameters, we fully consider the difficulty of measuring different types of data in humans, and adopt easily accessed data in non-invasive measurements, namely blood pressure and blood flow.
Moreover, given the mutual influence and connection between the various tissue structures of the human body, we explicitly model the cardiovascular system and the liver system through the blood vessel network, which can comprehensively simulate the spread of blood flow in the human body.

To the best of our knowledge, our study serves as the first to graft hemodynamic models with clinical measurements in human bodies. However, there are several limitations in this work.
First, from a medical point of view, although the hemodynamic model provides an explanation for the physiological behavior of the human body to a certain extent, the mathematical and physical models fail to fully consider the possible effects of various factors on blood transmission, such as the level of blood oxygen concentration and trace elements in the blood.
Therefore, the estimated parameters need to be adjusted appropriately, and the inclusive criteria of admitting volunteers are very strict, which limit the sample size of the study.
Second, the heart model in this article uses the pump model and the valve model to simulate the ventricles and valves, which actually greatly simplifies the complex heart system. In addition, this model is based on the heart data of normal people, and it has certain limitations for patients who have some problems with the cardiovascular system itself.

To provide reliable suggestions for clinical applications through modeling and numerical simulation, there are many aspects that are worthwhile for further research.
First, future work should consider the complex physiological behaviors of various organs and tissues in the human body, which are omitted in this work for the sake of simplicity.
Second, although the one-dimensional hemodynamic equation can simulate the blood flow propagation accurately, future research may consider the three-dimensional blood vessel model, which can fully reflect the influence of more factors on the blood propagation, to improve the estimation accuracy.
Third, due to the structural differences between the one-dimensional blood vessel model and the three-dimensional blood vessel model, the coupling results are not satisfactory; hence, we fail to use a one-dimensional blood vessel model to solve the boundary conditions of a three-dimensional blood vessel model. Therefore, it is worthwhile to examine three-dimensional model and the coupling methods, such as the terminal impedance model, so that this framework can adapt to more vascular network scenarios.
Furthermore, the blood vessel model can be extended to a more subtle level, for example, consider the establishment of a microscopic model of blood cell transmission to study the physiological characteristics of blood.
In order to meet the needs of clinical surgery, it is also necessary to study how to add various surgical behaviors such as incisions, clipping and other operations into the model to adapt it to a variety of surgical scenarios. In addition, the drugs and adjusted positions used in surgery  also have a certain impact on the operation of various organs of the human body; hence are of great value to be considered in future studies.



 \section*{Acknowledgments}
 We would like to acknowledge the assistance of volunteers in putting
 together this example manuscript and supplement.

\bibliographystyle{siamplain}
\bibliography{references}
\end{document}

%% file: ex_shared.tex
\usepackage{lipsum}
\usepackage{amsfonts}
\usepackage{graphicx}
\usepackage{epstopdf}
\usepackage{algorithmic}
\ifpdf
  \DeclareGraphicsExtensions{.eps,.pdf,.png,.jpg}
\else
  \DeclareGraphicsExtensions{.eps}
\fi

\newsiamremark{remark}{Remark}
\newsiamremark{hypothesis}{Hypothesis}
\crefname{hypothesis}{Hypothesis}{Hypotheses}
\newsiamthm{claim}{Claim}

\headers{Data Analysis On Human Liver Circulation}{T. Wu, E. Fainman, and Y. Wang}

\title{A Data Analysis Study On Human Liver Blood Circulation\thanks{Submitted to the editors DATE.
}}

\author{Ting Wu\thanks{Department of Mathematics, Nanjing University 
  (\email{tingwu@nju.edu.cn})}
\and  Keqin Liu\thanks{Department of Mathematics, Nanjing University 
  (\email{kqliu@nju.edu.cn})}
\and Emily Fainman\thanks{McCoy College of Business, Texas State University
  (\email{c\_z88@txstate.edu})}
\and Ying Wang\thanks{Department of Anesthesia, The Fourth Affiliated Hospital of Harbin Medical University}
}

\usepackage{amsopn}


%% file: main.bbl
\begin{thebibliography}{10}

\bibitem{Audebert2017partial}
{\sc C.~Audebert, M.~Bekheit, P.~Bucur, E.~Vibert, and I.~E. Vignon-Clementel},
  {\em Partial hepatectomy hemodynamics changes: Experimental data explained by
  closed-loop lumped modeling}, Journal of biomechanics, 50 (2017),
  pp.~202--208.

\bibitem{Audebert2017}
{\sc C.~Audebert, P.~Bucur, M.~Bekheit, E.~Vibert, I.~E. Vignon-Clementel, and
  J.-F. Gerbeau}, {\em Kinetic scheme for arterial and venous blood flow, and
  application to partial hepatectomy modeling}, Computer Methods in Applied
  Mechanics and Engineering, 314 (2017), pp.~102--125.

\bibitem{Audusse2000}
{\sc E.~Audusse, M.-O. Bristeau, and B.~Perthame}, {\em Kinetic schemes for
  Saint-Venant equations with source terms on unstructured grids}, PhD thesis,
  INRIA, 2000.

\bibitem{barrett2010ganong}
{\sc K.~E. Barrett, S.~Boitano, S.~M. Barman, and H.~L. Brooks}, {\em
  Ganong’s review of medical physiology twenty},  (2010).

\bibitem{Berzigotti2004}
{\sc A.~Berzigotti, S.~Dapporto, L.~Angeloni, S.~Ramilli, G.~Bianchi, M.~C.
  Morelli, D.~Magalotti, and M.~Zoli}, {\em Postprandial splanchnic
  haemodynamic changes in patients with liver cirrhosis and patent
  paraumbilical vein}, European journal of gastroenterology \& hepatology, 16
  (2004), pp.~1339--1345.

\bibitem{Boileau2015}
{\sc E.~Boileau, P.~Nithiarasu, P.~J. Blanco, L.~O. M{\"u}ller, F.~E. Fossan,
  L.~R. Hellevik, W.~P. Donders, W.~Huberts, M.~Willemet, and J.~Alastruey},
  {\em A benchmark study of numerical schemes for one-dimensional arterial
  blood flow modelling}, International journal for numerical methods in
  biomedical engineering, 31 (2015), p.~e02732.

\bibitem{Bonfiglio2010}
{\sc A.~Bonfiglio, K.~Leungchavaphongse, R.~Repetto, and J.~H. Siggers}, {\em
  Mathematical modeling of the circulation in the liver lobule}, Journal of
  biomechanical engineering, 132 (2010).

\bibitem{Brook2002}
{\sc B.~Brook and T.~Pedley}, {\em A model for time-dependent flow in (giraffe
  jugular) veins: uniform tube properties}, Journal of biomechanics, 35 (2002),
  pp.~95--107.

\bibitem{Childress2014}
{\sc E.~M. Childress and C.~Kleinstreuer}, {\em Impact of fluid--structure
  interaction on direct tumor-targeting in a representative hepatic artery
  system}, Annals of biomedical engineering, 42 (2014), pp.~461--474.

\bibitem{Chu1992}
{\sc T.-M. Chu and N.~P. Reddy}, {\em A lumped parameter mathematical model of
  the splanchnic circulation}, J Biomech Eng., 114 (1992), pp.~222--226.

\bibitem{Colciago2014}
{\sc C.~M. Colciago, S.~Deparis, and A.~Quarteroni}, {\em Comparisons between
  reduced order models and full 3d models for fluid--structure interaction
  problems in haemodynamics}, Journal of Computational and Applied Mathematics,
  265 (2014), pp.~120--138.

\bibitem{Euler1763}
{\sc L.~Euler}, {\em Principia theoriae machinarum}, Novi Commentarii academiae
  scientiarum Petropolitanae,  (1763), pp.~230--253.

\bibitem{Formaggia2003}
{\sc L.~Formaggia, D.~Lamponi, and A.~Quarteroni}, {\em One-dimensional models
  for blood flow in arteries}, Journal of engineering mathematics, 47 (2003),
  pp.~251--276.

\bibitem{Formaggia1999}
{\sc L.~Formaggia, F.~Nobile, A.~Quarteroni, and A.~Veneziani}, {\em Multiscale
  modelling of the circulatory system: a preliminary analysis}, Computing and
  visualization in science, 2 (1999), pp.~75--83.

\bibitem{Garcea2009}
{\sc G.~Garcea and G.~Maddern}, {\em Liver failure after major hepatic
  resection}, Journal of hepato-biliary-pancreatic surgery, 16 (2009),
  pp.~145--155.

\bibitem{Fitzgerald2018}
{\sc F.~Grace}, {\em The aorta}.
\newblock Teach Me Anatomy, 2018,
  \url{https://teachmeanatomy.info/abdomen/vasculature/arteries/aorta/}
  (accessed 2021-07-13).

\bibitem{Ho2012}
{\sc H.~Ho, K.~Sorrell, A.~Bartlett, and P.~Hunter}, {\em Blood flow simulation
  for the liver after a virtual right lobe hepatectomy}, in International
  Conference on Medical Image Computing and Computer-Assisted Intervention,
  Springer, 2012, pp.~525--532.

\bibitem{Janela2010}
{\sc J.~Janela, A.~Moura, and A.~Sequeira}, {\em Absorbing boundary conditions
  for a 3d non-newtonian fluid--structure interaction model for blood flow in
  arteries}, International Journal of Engineering Science, 48 (2010),
  pp.~1332--1349.

\bibitem{Lagana2002}
{\sc K.~Lagana, G.~Dubini, F.~Migliavacca, R.~Pietrabissa, G.~Pennati,
  A.~Veneziani, and A.~Quarteroni}, {\em Multiscale modelling as a tool to
  prescribe realistic boundary conditions for the study of surgical
  procedures}, Biorheology, 39 (2002), pp.~359--364.

\bibitem{Markou2004}
{\sc N.~Markou, L.~Grigorakos, P.~Myrianthefs, E.~Boutzouka, M.~Rizos,
  P.~Evagelopoulou, H.~Apostolakos, and G.~Baltopoulos}, {\em Venous pressure
  measurements in the superior and inferior vena cava: the influence of
  intra-abdominal pressure.}, Hepato-gastroenterology, 51 (2004), pp.~51--55.

\bibitem{martini2006human}
{\sc F.~Martini, M.~J. Timmons, R.~B. Tallitsch, W.~C. Ober, C.~W. Garrison,
  K.~B. Welch, and R.~T. Hutchings}, {\em Human anatomy}, Pearson/Benjamin
  Cummings San Francisco, CA, 2006.

\bibitem{Papadakis2009}
{\sc G.~Papadakis}, {\em Coupling 3d and 1d fluid--structure-interaction models
  for wave propagation in flexible vessels using a finite volume
  pressure-correction scheme}, Communications in numerical methods in
  engineering, 25 (2009), pp.~533--551.

\bibitem{Pedley1996}
{\sc T.~J. Pedley, B.~S. Brook, and R.~S. Seymour}, {\em Blood pressure and
  flow rate in the giraffe jugular vein}, Philosophical Transactions of the
  Royal Society of London. Series B: Biological Sciences, 351 (1996),
  pp.~855--866.

\bibitem{Peeters2015}
{\sc G.~Peeters, C.~Debbaut, P.~Cornillie, T.~De~Schryver, D.~Monbaliu,
  W.~Laleman, and P.~Segers}, {\em A multilevel modeling framework to study
  hepatic perfusion characteristics in case of liver cirrhosis}, Journal of
  biomechanical engineering, 137 (2015).

\bibitem{Perthame2002}
{\sc B.~Perthame et~al.}, {\em Kinetic formulation of conservation laws},
  vol.~21, Oxford University Press, 2002.

\bibitem{Quarteroni2003}
{\sc A.~Quarteroni and A.~Veneziani}, {\em Analysis of a geometrical multiscale
  model based on the coupling of ode and pde for blood flow simulations},
  Multiscale Modeling \& Simulation, 1 (2003), pp.~173--195.

\bibitem{Ricken2015}
{\sc T.~Ricken, D.~Werner, H.~Holzh{\"u}tter, M.~K{\"o}nig, U.~Dahmen, and
  O.~Dirsch}, {\em Modeling function--perfusion behavior in liver lobules
  including tissue, blood, glucose, lactate and glycogen by use of a coupled
  two-scale pde--ode approach}, Biomechanics and modeling in mechanobiology, 14
  (2015), pp.~515--536.

\bibitem{Riemann1892}
{\sc B.~Riemann, R.~Dedekind, and H.~Weber}, {\em Gesammelte mathematische
  Werke und wissenschaftlicher Nachlass}, Рипол Классик, 1892.

\bibitem{rubin1998measurement}
{\sc G.~D. Rubin, D.~S. Paik, P.~Johnston, and S.~Napel}, {\em Measurement of
  the aorta and its branches with helical ct.}, Radiology, 206 (1998),
  pp.~823--829.

\bibitem{Schaaf1972}
{\sc B.~W. Schaaf and P.~H. Abbrecht}, {\em Digital computer simulation of
  human systemic arterial pulse wave transmission: a nonlinear model}, Journal
  of Biomechanics, 5 (1972), pp.~345--364.

\bibitem{Sherwin2003}
{\sc S.~Sherwin, V.~Franke, J.~Peir{\'o}, and K.~Parker}, {\em One-dimensional
  modelling of a vascular network in space-time variables}, Journal of
  engineering mathematics, 47 (2003), pp.~217--250.

\bibitem{Sherwin2003a}
{\sc S.~J. Sherwin, L.~Formaggia, J.~Peiro, and V.~Franke}, {\em Computational
  modelling of 1d blood flow with variable mechanical properties and its
  application to the simulation of wave propagation in the human arterial
  system}, International journal for numerical methods in fluids, 43 (2003),
  pp.~673--700.

\bibitem{Stettler1981}
{\sc J.~Stettler, P.~Niederer, and M.~Anliker}, {\em Theoretical analysis of
  arterial hemodynamics including the influence of bifurcations}, Annals of
  biomedical engineering, 9 (1981), pp.~145--164.

\bibitem{Suzuki2019}
{\sc T.~Suzuki, Y.~Suzuki, J.~Okuda, R.~Minoshima, Y.~Misonoo, T.~Ueda,
  J.~Kato, H.~Nagata, T.~Yamada, and H.~Morisaki}, {\em Cardiac output and
  stroke volume variation measured by the pulse wave transit time method: a
  comparison with an arterial pressure-based cardiac output system}, Journal of
  clinical monitoring and computing, 33 (2019), pp.~385--392.

\bibitem{Toro2013}
{\sc E.~F. Toro and A.~Siviglia}, {\em Flow in collapsible tubes with
  discontinuous mechanical properties: mathematical model and exact solutions},
  Communications in Computational Physics, 13 (2013), pp.~361--385.

\bibitem{VanDerPlaats2004}
{\sc A.~Van Der~Plaats, N.~'tHart, G.~Verkerke, H.~Leuvenink, P.~Verdonck,
  R.~Ploeg, and G.~Rakhorst}, {\em Numerical simulation of the hepatic
  circulation}, The International journal of artificial organs, 27 (2004),
  pp.~222--230.

\bibitem{Westerhof2009}
{\sc N.~Westerhof, J.-W. Lankhaar, and B.~E. Westerhof}, {\em The arterial
  windkessel}, Medical \& biological engineering \& computing, 47 (2009),
  pp.~131--141.

\bibitem{Womersley1957}
{\sc J.~R. Womersley}, {\em An elastic tube theory of pulse transmission and
  oscillatory flow in mammalian arteries}, tech. report, Aerospace Research
  Labs Wright-Patterson AFB OH, 1957.

\end{thebibliography}
